\documentclass[12pt,a4paper,absolute]{mymathart}
\usepackage{mymathmacros}

\newcommand{\Ek}{E_{\kappa}}

\include{hyphenation}

\usepackage{graphicx}
\usepackage{subfigure}

\usepackage{hyperref}

\theoremstyle{changebreak}

\newtheorem{thm}{Theorem}[section]
\newtheorem{defn}[thm]{Definition}
\newtheorem{lem}[thm]{Lemma}
\newtheorem{cor}[thm]{Corollary}

\newtheorem{obs}[thm]{Observation}

\swapnumbers

\newtheorem{mainthm}{Theorem}
\newtheorem*{maincor}{Corollary}
\newtheorem{appthm}[thm]{Theorem}
\newtheorem{applem}[thm]{Lemma}
\newtheorem{appcor}[thm]{Corollary}

\newcommand{\W}{\mathcal{W}}

\newcommand{\T}{\mathcal{T}}
\newcommand{\K}{\mathcal{K}}

\renewcommand{\P}{\mathcal{P}}
\newcommand{\B}{\mathcal{B}}

\newcommand{\Herman}{\mathcal{H}}

\renewcommand{\phi}{\varphi}

\newcommand{\id}{\operatorname{id}}

\title[Siegel Disks and Periodic Rays]{%
  Siegel Disks and Periodic Rays\\ of Entire Functions}

\author{Lasse Rempe}
\address{Department of Mathematical Sciences, University of Liverpool,
L69 7ZL, UK}
\email{l.rempe@liverpool.ac.uk}
\date{today}
\thanks{Supported by a postdoctoral fellowship of the 
 German Academic Exchange Service (DAAD) and
 by the German-Israeli Foundation
 for Scientific Research and Development (G.I.F.),
 grant no.\ G-643-117.6/1999}

\subjclass{Primary 37F10; Secondary 30D05}

\date{\today}

\begin{document}

 \begin{abstract}
  Let $f$ be an entire function whose set of singular values is bounded
  and suppose that $f$ has a Siegel disk $U$ such that 
  $f|_{\partial U}$ is a homeomorphism. We show that $U$ is bounded.
  Using a result of Herman, we deduce that if additionally 
  the rotation number of $U$ is Diophantine, then $\partial U$ contains
  a critical point of $f$. 

  Suppose furthermore that all singular values of $f$ lie in the Julia set.
   We prove that, if $f$ has a Siegel disk $U$ whose boundary
   contains no singular values, then the condition
   that $f\colon\partial U\to\partial U$ 
   is a homeomorphism is automatically satisfied. 
   We also 
   investigate landing properties of periodic dynamic rays
   by similar methods.
 \end{abstract}

\maketitle

 \section{Introduction}

 \subsection*{Main Results}

  Let $f\colon\C\to\C$ be a nonlinear entire function and suppose that
   $U$ is a Siegel disk of $f$, i.e.\ an invariant component of the
   Fatou set on which $f$ is conformally conjugate to an irrational
   rotation. It is an important question under which conditions
   $\partial U$ contains a (finite) singular value of $f$. 

  Herman \cite{hermansiegel} proved the following result. 
    \emph{If the rotation number of $U$ is Diophantine\footnote{%
    Diophantine rotation numbers can be replaced by the larger
    class $\Herman$ described by Yoccoz \cite{yoccozcircle};
    compare the discussion in \cite[Chapter I]{perezmarco}.},
    if $U$ is bounded
    and if $f\colon{\partial U}\to\partial U$ is a homeomorphism, 
    then $\partial U$ contains a 
    critical point.}

\begin{figure}
 \subfigure[$f(z)=\lambda(\exp(z)-1)$]{%
   \fbox{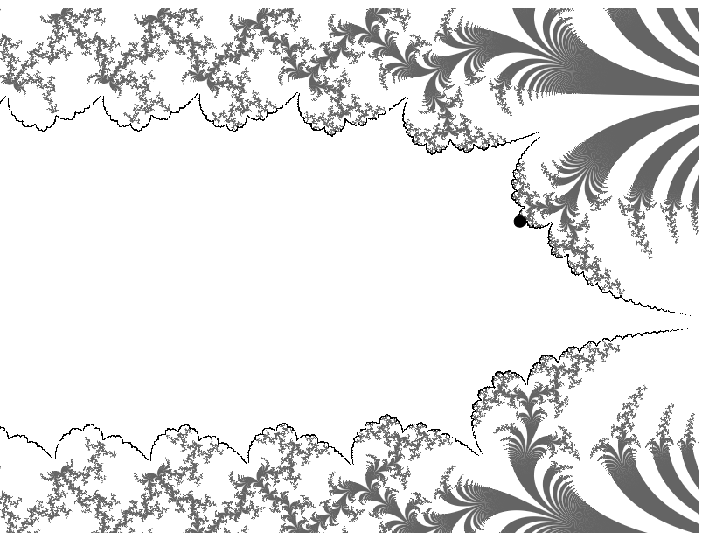}\label{fig:lambdaexp}}\hfill%
 \subfigure[$f(z)=\lambda\sin(z)$]{%
   \fbox{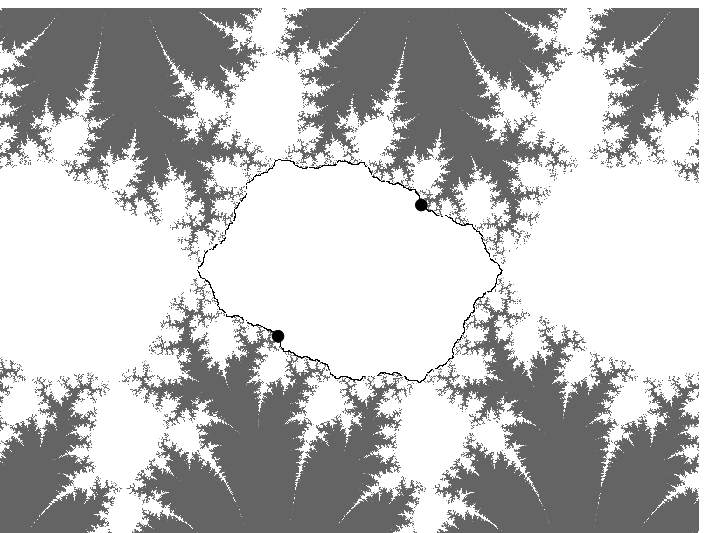}\label{fig:lambdasin}}
 \caption{Two entire functions with a Diophantine Siegel disk
   (in both cases the rotation number is the golden mean). Our results
   show that both Siegel disk boundaries do indeed contain a
   singular value.}
\end{figure}

   Here the condition that $f|_{\partial U}$ be a homeomorphism
    is needed to exclude certain topological ``pathologies''
     (compare
    \cite{rogers});
    it is currently unknown whether these can actually occur.
    With this reasonable restriction, Herman's theorem gives
    a very satisfactory answer when $f$ is a polynomial, since
    $U$ is always bounded in this case. On the other hand,
    the result is far less complete for transcendental functions.
    For example,
    it does not answer the question whether the boundary of the
    Siegel disk depicted in Figure \ref{fig:lambdaexp} does indeed contain
    the singular value. In this article, we prove 
    that the assumption of boundedness can be removed
    for a large class of
    entire functions.

   \begin{mainthm}[Univalent Siegel disks] \label{mainthm:univalentsiegel}
    Let $f\colon\C\to\C$ be an entire function whose set of
     singular values is bounded and suppose
     that $U$ is a Siegel disk of $f$ 
     with the property
     that $f\colon\partial U\to\partial U$ is a homeomorphism.
     Then $U$ is bounded. 

    In particular, if the rotation number of
     $U$ is Diophantine, then $\partial U$ contains a critical point. 
   \end{mainthm}
   \begin{remark}
    The class of entire functions whose set of singular values is bounded 
     is called the \emph{Eremenko-Lyubich class} and 
     is commonly denoted by $\B$. This 
     is the natural class of entire functions to which our methods
     apply. It seems by no means clear whether one
     should expect Siegel disk boundaries to contain finite singular
     values for more general classes of maps.
   \end{remark}

  Since the condition of Theorem \ref{mainthm:univalentsiegel} is often
   difficult to verify directly, we also prove the following.  
    (Here $J(f)$ and $S(f)$ denote the Julia set and set of singular
   values of $f$, respectively; compare the remarks on notation below.)
  \begin{mainthm}[Nonsingular Siegel disks]
   \label{mainthm:nonsingularsiegel}
   Let $f\in\B$ with $S(f)\subset J(f)$. If $U$ is a Siegel disk of $f$
   which satisfies $S(f)\cap \partial U = \emptyset$, then
   $f\colon\partial U\to \partial U$ is a homeomorphism.
  \end{mainthm}
   \begin{remark}
    If $f$ is an exponential map, i.e. $f(z) = e^{2\pi i \theta}(\exp(z)-1)$,
     and $f$ has a Siegel disk, then the unique singular value of $f$ 
     must automatically belong to the Julia set.
     So Theorems \ref{mainthm:univalentsiegel} and 
     \ref{mainthm:nonsingularsiegel} imply that the boundary of
     an unbounded Siegel disk
     of an exponential map always contains the singular value, which
     answers a question of Herman, Baker and
     Rippon \cite[Problem 2.86 (b)]{haymanlist}. The proof for
     this special
     case has previously appeared in \cite{siegel} and was 
     obtained independently by Buff and Fagella \cite{bufffagella}.
   \end{remark}

  Combined with Herman's result mentioned above, the preceding theorems 
   yield the following corollary. 
  \begin{maincor}[Diophantine Siegel disks]
  Let $f\in\B$ with $S(f)\subset J(f)$, and suppose that $f$ has a Siegel disk
   $U$ with diophantine rotation number. 
   Then $S(f)\cap \partial U\neq\emptyset$.
  \end{maincor}

  Finally, if $f$ has only two critical values and no asymptotic values,
   we can give a complete result on the boundedness of
   nonsingular Siegel disks. (An important example of such functions
    is given by the family
    $C_{a,b}\colon z\mapsto a\exp(z)+b\exp(-z)$ of cosine maps, 
    where $a,b\in \C\setminus\{0\}$.)

  \begin{mainthm}[Maps with two critical values]
   \label{mainthm:nonsingularcosinesiegel}
   Let $f$ be an entire function which
    has two critical values and no asymptotic values,
    and suppose that $f$ has a 
   periodic (i.e., not necessarily 
   invariant) Siegel disk $U$ such that, for all 
   $j\geq 0$, the boundary $\partial f^j(U)$ contains no 
   critical values of $f$. 
   Then $U$ is bounded.
  \end{mainthm}
 \begin{remark}[Remark 1]
  In the case where both critical values of $f$ lie in the
   Julia set, it is sufficient to demand that no $\partial f^j(U)$ contains
   both of these.
 \end{remark}
 \begin{remark}[Remark 2]
  It does not seem unreasonable
   to expect that all Siegel disks in the cosine family are bounded,
   but the difficulties involved in showing this
   are unresolved even for
   polynomials. For instance, it is not known whether a quadratic polynomial
   can have a Siegel disk whose boundary is the entire Julia set. A cosine
   Siegel disk whose boundary is the Julia set would in particular be
   unbounded. 
 \end{remark}

\smallskip
  The methods which yield the above theorems apply, in fact, to a wide
  range of connected invariant sets. As a second example,
  we discuss the application to periodic rays. 
  A \emph{fixed ray} of an entire function $f$ is a curve  
    \[ \gamma\colon(-\infty,\infty)\to \C \]
  with $\lim_{t\to+\infty}|\gamma(t)|=\infty$ 
  which satisfies $f(\gamma(t))=\gamma(t+1)$ for 
  all $t$.
  As usual, we say that $\gamma$ \emph{lands} at a point
  $z_0\in\Ch$ if $\lim_{t\to -\infty}\gamma(t)=z_0$. 
  A \emph{periodic ray} of $f$ is a fixed ray of some iterate
   $f^n$ of $f$.

  Periodic rays play an essential part in the study of
   polynomial dynamics. It is now known that such rays exist
   for many, although not for all,
   functions $f$ in the Eremenko-Lyubich class,
   in particular for those of finite order \cite{fatoueremenko}.
   They have already been used
   to great advantage in the theory of exponential maps,
   see e.g.\ \cite{expper,expattracting}, and are likely to be equally useful
   in the study of more general families. 

  Much of the usefulness of periodic rays for polynomials
   depends on the fact that every such ray lands at a 
   repelling or parabolic
   periodic point  (see \cite[Theorem 18.10]{jackdynamicsthird} 
    or our Corollary
    \ref{cor:polylanding}). 
   While the proof of this result breaks down in the
   transcendental case (compare Appendix \ref{app:landing1}), the author
   was
   recently able to generalize it to exponential maps
   \cite{landing2new}. This gives us reason to 
   expect that such a result remains true
   for some larger classes of entire functions.
   However, the proof 
   in \cite{landing2new}
   uses a theorem of
   Schleicher \cite{habil} on exponential parameter space which depends
   essentially on the fact that this parameter space is one-dimensional. 
   We will prove the following result.

  \begin{mainthm}[Landing of periodic rays]
   \label{mainthm:landing}
   Let $f\in\B$. If $\gamma$ is a fixed ray of $f$ such that
    $f\colon\cl{\gamma}\to\cl{\gamma}$ 
    is a homeomorphism and the accumulation set of
    $\gamma$ does not contain any critical points, then
    $\gamma$ lands at a repelling or parabolic fixed point of $f$.

   Similarly, if $S(f)\subset J(f)$ and $\gamma$ is a fixed ray  
    of $f$ whose accumulation set does not intersect $S(f)$, then
   $\gamma$ lands at a repelling or parabolic fixed point of $f$.
  \end{mainthm}
  As far as we know, this is the first landing
   criterion not relying on hyperbolic expansion which can
   be applied to functions
   in higher-dimensional
   parameter spaces.

\subsection*{Idea and Structure of the Proof}

 As already mentioned, our results are not restricted to 
  the special (though important) cases described
  above. Their basis lies in the following general principle.

 \begin{mainthm}[Invariant connected sets] \label{mainthm:pullbacks}
  Let $f\in\B$ and suppose that $A\subset\C$ is closed and connected
   such that
   $f(A)\subset A$ and such that
   $f\colon A\to \cl{f(A)}$ is a homeomorphism.

  Then for every $R>0$, there exists $R'>0$ such that
\[
     \{z\in A\colon |z|\geq R'\} \subset
       \{z\in \C\colon 
         f^n(z)\geq R \text{ for all $n$}
        \text{and }|f^n(z)|\to\infty \text{ as $n\to\infty$}\}. \]
   \end{mainthm}

  This result implies Theorem \ref{mainthm:univalentsiegel}
   by letting $A$ be the closure of the Siegel disk $U$ and
   using the fact that 
   $U$ cannot contain escaping points. (See
   Section \ref{sec:mainpullbacks} for details.)    

  The proof of Theorem \ref{mainthm:pullbacks} is,
   in fact, quite simple. The hypothesis implies that 
   the unbounded parts of $A$ are contained in 
   finitely many \emph{fundamental domains} of $f$ (for cosine
   maps $C_{a,b}$ as above
   this is equivalent to $\im A$ being bounded). An expansion 
   argument then shows that any sufficiently large point in
   $A$ has a large image (this is akin to the fact that
   $|C_{a,b}(z)|$ behaves like 
   $\exp(|\re z|)$ when $|\re z|$ is large), yielding the desired result. 

  This proof is carried out in Sections
   \ref{sec:expansion} and \ref{sec:mainpullbacks}, with the former
   section
   reviewing basic definitions for functions in the Eremenko-Lyubich class
   and deducing
   the abovementioned expansion statement, and the latter containing
   the actual proof. 

  The remaining two sections show how
   Theorem \ref{mainthm:pullbacks} can be applied in cases
   where there are no singular values in the Fatou set, and how to
   apply our results to the landing problem for periodic rays. 

  Two auxiliary results of a topological nature were relegated to Appendix
   \ref{app:topology} to avoid interrupting the flow of ideas. 
   Appendix \ref{app:landing1} discusses difficulties in proving
   landing results for periodic dynamic rays using hyperbolic
   contraction. 

\subsection*{Acknowledgments}
  I would like to thank Walter Bergweiler, Adam Epstein, Markus
   F\"orster, Jan Kiwi,
   Phil Rippon,
   G\"unter Rottenfu{\ss}er, Gwyneth Stallard and Sebastian van Strien
   as well as the audiences at the
   IHP, the Universities of Warwick and Liverpool and at the
   Max-Planck-Institut Bonn for interesting
   discussions. I would also like to thank the University of Warwick
   for its hospitality during the time this work was conducted. Finally,
   I would like to thank the referee for many helpful and detailed comments.

\subsection*{Notation}
  We denote the complex plane and the Riemann sphere
   by $\C$ and $\Ch := \C\cup\{\infty\}$, respectively. 
   The unit disk is denoted by $\D$; more generally, $\D_r(z_0)$ denotes 
   the open disk of radius $r$ around some $z_0\in\C$. 
   The unit circle is denoted by $S^1$. 

   The underlying topological space for our considerations is the 
    complex plane $\C$; all closures, boundaries, neighborhoods
    etc.\ will be understood
    to be taken in $\C$ unless explicitly stated otherwise. 
    The Euclidean length of a curve $\gamma$ is denoted by
    $\ell(\gamma)$

  Throughout this article (with the exception of 
   the first half of Section \ref{sec:landing}),
   $f\colon\C\to\C$ will be a nonconstant nonlinear entire function. As usual,
   the Fatou and Julia sets of $f$ are denoted by
   $F(f)$ and $J(f)$; the \emph{set of escaping points} of $f$ is 
    \[ I(f) := \{z\in\C\colon |f^n(z)|\to\infty\}. \] 

   The set of \emph{singularities of $f^{-1}$}, denoted
   $\sing(f^{-1})$, consists
   of all finite critical and asymptotic values of $f$; the elements of
    $S(f) := \cl{\sing(f^{-1})}$
   are called the \emph{singular values} of $f$. We will mostly be
   interested in the aforementioned \emph{Eremenko-Lyubich class}
    \[ \B := \{f\colon\C\to\C\text{ transcendental entire}\colon
                  S(f) \text{ is bounded}\}. \]

   An \emph{(invariant) Siegel disk} of $f$ is a simply connected
    component of $F(f)$
    with $f(U)=U$ such that $f|_U$ is conjugate to an irrational
    rotation. A \emph{periodic Siegel disk} is a component of
    $F(f)$ which is an invariant Siegel disk for some iterate of $f$.

   We conclude any proof by the symbol $\proofsymbol$. Proofs of
    separate claims within a larger proof will be completed by
    $\subproofsymbol$, while statements cited without proof are 
    indicated by $\noproofsymbol$.


\section{Tracts and Expansion}
   \label{sec:expansion}

 Throughout this section, we 
  fix some entire function $f$ in the Eremenko-Lyubich class $\B$.
  We will review some of the standard
  constructions used when dealing with such maps 
  and deduce an expansion property
  (Lemma \ref{lem:expansion} below) which is essential
  for our arguments. We also recall some recent results
  on the existence of unbounded connected
  subsets of $I(f)$.

 \subsection*{Tracts and Fundamental Domains}
  Define $K:=1+\max\bigl(|f(0)|,\max_{s\in S(f)}|s|\bigr)$ and  
  $G_K:=\{z\in\C\colon |z|>K\}$. Then each component of
  $f^{-1}(G_K)$ is a simply
  connected domain whose boundary is a Jordan arc tending to infinity in
  both directions.
  These components are called the \emph{tracts} of $f$; restricted
  to any tract, $f$ is a universal covering
  onto
  $G_K$. 

 Let $\gamma$ be a curve
  in $G_K$ which does not intersect any tracts and which connects
  $\partial G_K$ to $\infty$. (For example, we can let $\gamma$ be 
  a piece of the boundary of one of the tracts, together with a curve
  connecting it to $\partial G_K$ if necessary.)
  Then $f^{-1}(\gamma)$ cuts every tract
  into countably many components, which we call \emph{fundamental
  domains}; each fundamental domain maps univalently to
  $G_K\setminus\gamma$ under $f$. We would like to note the
  elementary but important fact that any bounded subset of $\C$
  intersects at most finitely many fundamental domains of $f$. 

 Our goal in this section is to prove the following lemma. 
  It basically states that a point whose orbit stays within finitely
  many fundamental domains escapes to infinity, provided it is large enough. 

  \begin{lem}[Growth of orbits]
    \label{lem:pointsgrowinfundamentaldomains}
    Suppose $F_1,\dots,F_k$ are fundamental domains of $f$, and
     $R>0$. Denote by $X$ the set of all $z\in\C$ with the following
     property:
   \begin{center}
    If $n\geq 0$ with $|f^m(z)|\geq R$ for $m=0,\dots,n$, then
     $\displaystyle{f^n(z)\in\bigcup_{j=1}^{k} F_j}$.
   \end{center}
     Then there exists $R'>0$ such that
      \[ X\cap \{|z|\geq R'\} \subset
         \{z\in I(f)\colon |f^n(z)|\geq R\text{ for all $n\geq 0$}\}. \]
  \end{lem}
  \begin{remark}
   $X$ consists of all points which do not leave the union of $F_1,\dots,F_k$
    without passing through the disk $\D_R(0)$ first; in particular, 
    $X$ contains this disk $\D_R(0)$ itself. $X$ will generally not be
    invariant under $f$, but if $z\in X$ and
    $|z|\geq R$, then $f(z)\in X$. 
  \end{remark}

 We wish to note that the subsequent constructions and results,
  which lead to the proof of this lemma and
  are perhaps 
  of a somewhat technical nature, are not required for the other
  sections of the article.

\begin{figure}
  \subfigure{%
\fbox{\resizebox{!}{.18\textheight}{%
   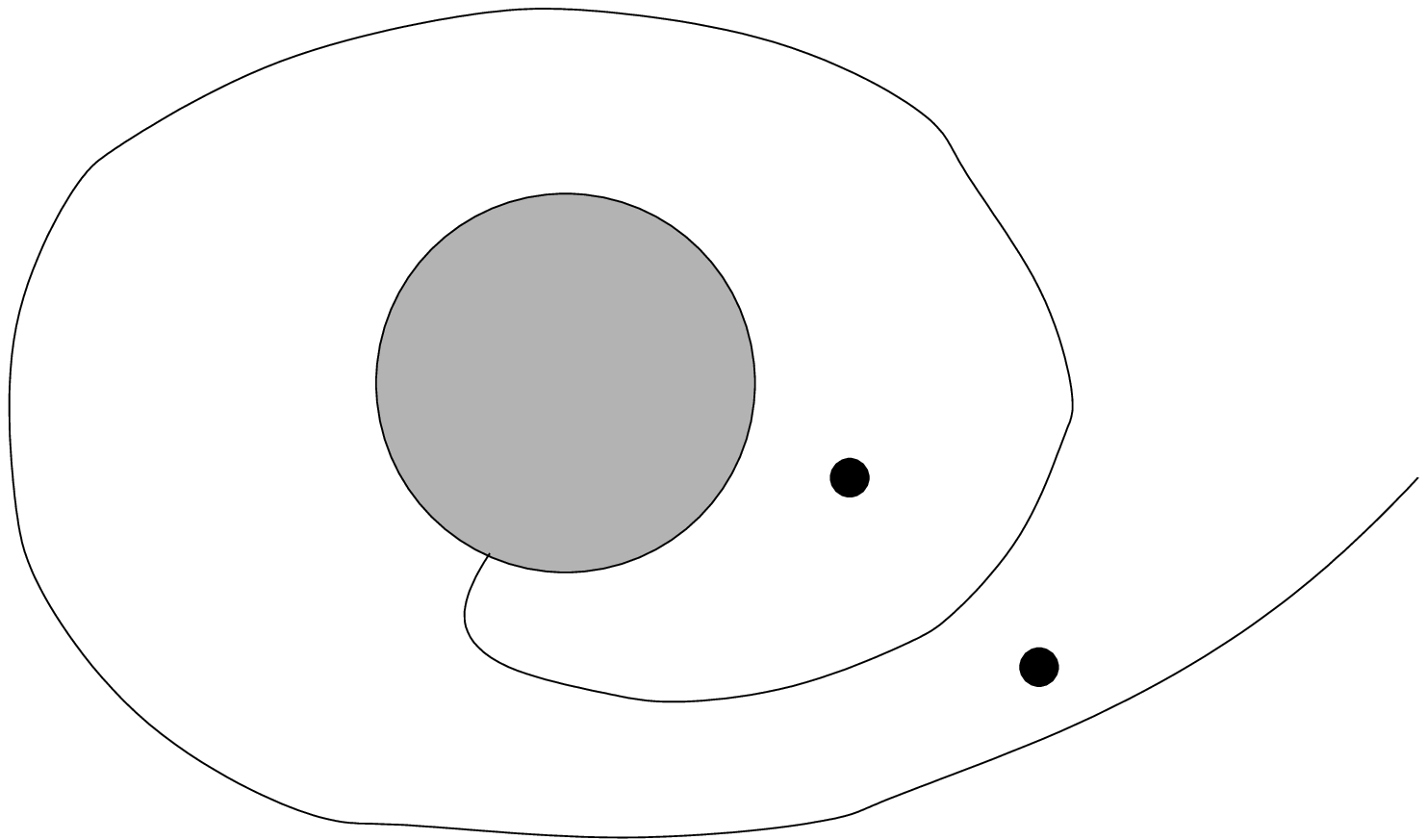}}}\hfill
  \subfigure{%
\fbox{\resizebox{!}{.18\textheight}{%
   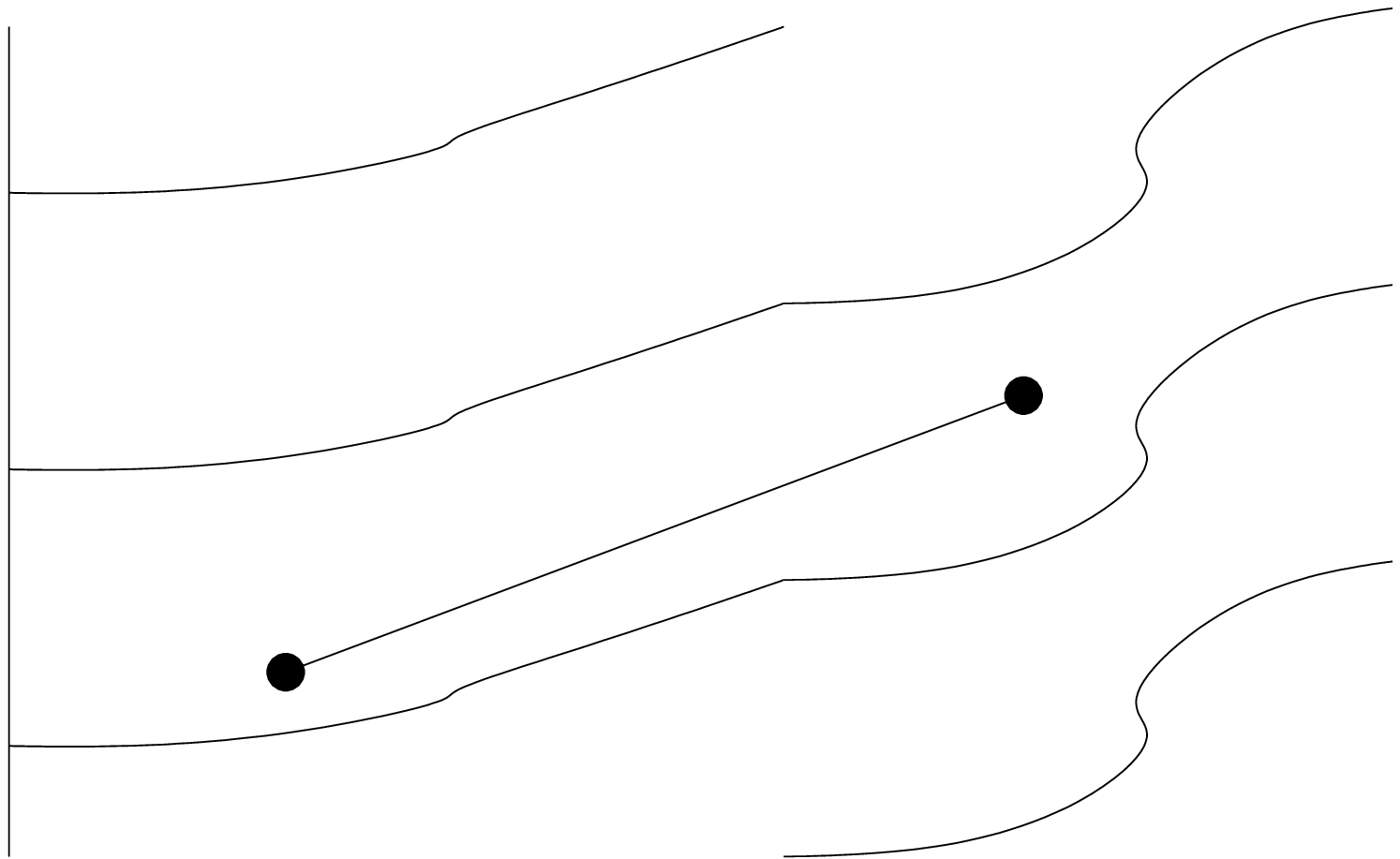}}}
\caption{Definition of $r(z) = \wt{r}(\zeta_0)$. On the left, the picture
  is in the $z$-plane, while the right hand side is drawn
  in logarithmic coordinates.\label{fig:definitionofr}} 
\end{figure}

 \subsection*{Logarithmic Coordinates}
  In \cite{alexmisha}, functions in class $\B$ were studied by applying
   a logarithmic change of variable both on $G_K$ and on the
   tracts of $f$. 
   More precisely, let $H := \{\re z > \log K\} = \exp^{-1}(G_K)$.
   Since $0\notin f^{-1}(G_K)$, and since $f$ is a universal covering
   on every tract, we can find a map $\Phi$ from the set
   $\T := \exp^{-1}(f^{-1}(G_K))$ to $H$ with 
   $\exp\circ\,\Phi = f\circ\exp$. If  $T$ is a component 
   of $\T$, then $T$ is simply connected and $\exp(T)$ is a tract of $f$;
   we call $T$ a tract of $\Phi$. The map
   $\Phi\colon T\to H$ 
    is a conformal isomorphism for every tract $T$ of $\Phi$. 

  Note that $\exp^{-1}(\gamma)$ consists of countably many curves,
  which cut the half plane $H$ into countably many \emph{fundamental
  strips}. By definition, the boundaries of these strips do not
  intersect any tracts of $\Phi$. 

  The following lemma --- proved by a simple application of Koebe's 
   $\frac{1}{4}$-theorem
   --- provides a basic expansion estimate for Eremenko-Lyubich functions.
   
  \begin{lem}[{\cite[Lemma 1]{alexmisha}}]\label{lem:alexmishaexpansion}
   For any $\zeta\in \T$, 
    \[
       |\Phi'(\zeta)| \geq \frac{1}{4\pi} (\re\Phi(\zeta) -\log K).
       \qedd \]
  \end{lem}

\subsection*{Growth of points in a fundamental domain} 
  In order to prove Lemma \ref{lem:pointsgrowinfundamentaldomains},
   we wish to show 
   that, for any fixed fundamental domain,
   every ``sufficiently large'' point has an even larger image. 
   While we
   cannot expect this statement to hold when size is measured
   by the modulus of a point, we will associate a 
   size $r(z)$ to a point $z$ which makes it true. 

  In order to do this, let us fix any base point
     $z_0\in G_{e^{16\pi}K}\setminus \gamma$
    for the remainder of this section. 
    If 
    $z\in G_K\setminus\gamma$ and
    $\zeta \in \exp^{-1}(z)$, we define
   $r(z) := \wt{r}(\zeta) := |\zeta - \zeta_0|$,
  where $\zeta_0$ is the unique point of $\exp^{-1}(z_0)$ which belongs
  to the same fundamental strip as $\zeta$. (Compare Figure
    \ref{fig:definitionofr}.)
  Note that
  $r(z)\geq \log|z| - \log|z_0|$, and that $r(z)$ remains bounded
  when $z$ ranges over a bounded subset of $G_K\setminus\gamma$.
  (However, $r(z)$ may tend to infinity 
   much faster than $\log|z|$).

\begin{figure}
 \resizebox{.99\textwidth}{!}{%
   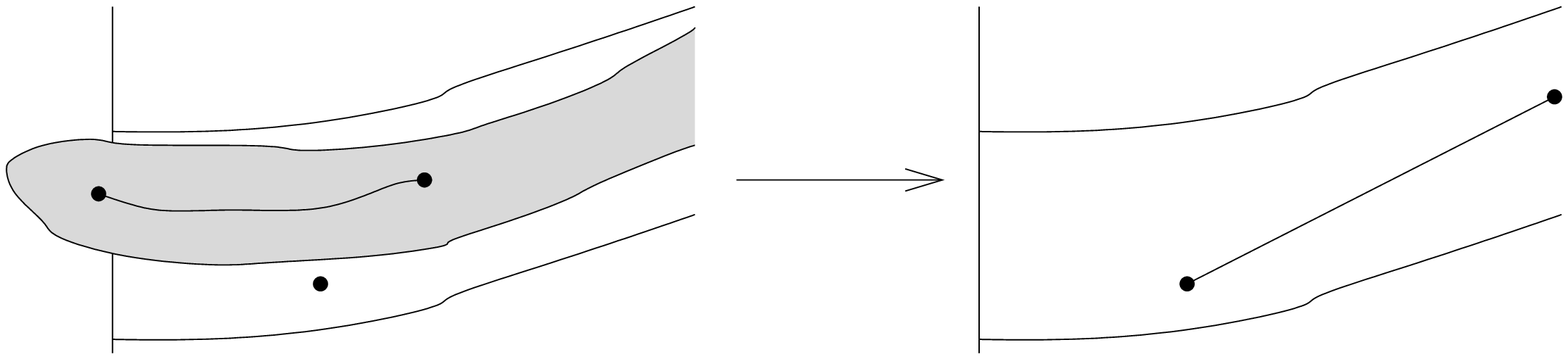}
\caption{Proof of Lemma \ref{lem:expansion}.\label{fig:proofoflemma}} 
\end{figure}

 \begin{lem}[Expansion in fundamental domains] \label{lem:expansion}
  Let $F$ be any fundamental domain. Then there exists some $C_F>0$
  with the following property. If $z\in F\cap G_K$ with $r(z)\geq C_F$,
  then $r(f(z))\geq 2r(z)$.
 \end{lem}
 \begin{proof} The proof will take place completely in logarithmic
   coordinates. Let $S$ be any fundamental strip and let
   $\zeta_0$ be the unique point of $\exp^{-1}(z_0)\cap S$. 
   Let $\wt{F}$ be the unique component
    of $\exp^{-1}(F)$ for which
    $S\cap \wt{F}$ is unbounded, and let $\wt{T}$
    be the tract of $\Phi$ containing $\wt{F}$. Then $\Phi$ maps
    $\wt{F}$ conformally onto some fundamental strip; by postcomposing
    $\Phi$ with a suitable translation we may suppose, for simplicity,
    that $\Phi(\wt{F})=S$. For the remainder of the proof,
    we denote the inverse $(\Phi|_{\wt{T}})^{-1}$ simply by $\Phi^{-1}$.

    Set $\omega_0 := \Phi^{-1}(\zeta_0)\in \wt{T}$ and define
   \begin{align*}
      C_1 &:= 1+\sup\bigl\{\wt{r}(\zeta):
                \zeta\in H\cap\wt{F} \text{ and }            
                \re\Phi(\zeta) < \log K + 16\pi\bigr\},
      \\
      C_2 &:= 1+\sup\bigl\{\wt{r}(\zeta): 
                            \zeta\in (H \cap \wt{F}) \setminus S\bigr\}
    \end{align*}
   and $C_F := \max( C_1, C_2, 2\cdot|\omega_0 - \zeta_0|)$. 
   We must show that
    \[ |\Phi(\zeta) - \zeta_0| \geq 2 \wt{r}(\zeta) \]
   for every point $\zeta\in \wt{F}\cap H$ with $\wt{r}(\zeta) \geq C_F$.

  So suppose that $\zeta$ is such a point. 
   Then, by definition of $C_1$ and $C_2$,
   $\zeta\in S$
   and $\re \Phi(\zeta) \geq \log K + 16\pi$. 
   Let $\alpha$ denote the straight line segment connecting
   $\Phi(\zeta)$ and $\zeta_0$, and set $\beta := \Phi^{-1}(\alpha)$. 
   (See Figure \ref{fig:proofoflemma}.) Then,
   by Lemma \ref{lem:alexmishaexpansion}, $|\Phi'(z)|\geq 4$ for every
   point $z\in\beta$, and thus
   \[ \ell(\beta) \leq \frac{1}{4} \ell(\alpha) = 
          \frac{|\Phi(\zeta) - \zeta_0|}{4}. \]
   Since $\beta$ is a curve connecting $\zeta$ and $\omega_0$, it follows that
    \[ \wt{r}(\zeta) = |\zeta - \zeta_0| \leq
       |\zeta - \omega_0| + |\omega_0 - \zeta_0| \leq
       \frac{|\Phi(\zeta) - \zeta_0|}{4} + \frac{\wt{r}(\zeta)}{2}, \]
   which means that $|\Phi(\zeta) - \zeta_0|\geq 2\wt{r}(\zeta)$, as 
   required. \end{proof}

\begin{proof}[Proof of Lemma \ref{lem:pointsgrowinfundamentaldomains}]
  By Lemma \ref{lem:expansion}, for each $j$ there exists some
   $C_j>0$ such that $r(f(z))\geq 2r(z)$ for every 
   $z\in F_j\cap G_K$ with $r(z)\geq C_j$.
   Choose a $T>0$ such that 
    $|z|\geq R$ whenever $r(z)\geq T$ 
    and set 
    $C := \max\bigl(T,\max_{j} C_j\bigr)$.

   If $z\in X$ 
    with $r(z)\geq C$, then $f(z)\in X$ and $r(f(z))\geq 2r(z)$. 
    Indeed, we have $r(z)\geq T$, and hence $|z|\geq R$. By definition of
    $X$, this implies that $f(z)\in X$, and that $z\in F_j$ for some $j$.
    By the definition of $C$, we have $r(f(z))\geq 2r(z)$ as claimed. 

   It follows inductively that $|f^n(z)| \geq R$ and 
     $r(f^n(z)) \geq 2^n r(z)$ for all $n$; 
    in particular, $z\in I(f)$. The claim follows by choosing $R'\geq K$
    sufficiently large so that $r(z)\geq C$ whenever $|z|\geq R'$. 
\end{proof}

\subsection*{Continua Consisting of Escaping Points} 
  It is an open question, posed by Eremenko \cite{alexescaping}, whether,
  for every transcendental entire function,
  every component of $I(f)$ is unbounded. Recently, Rippon and
  Stallard \cite{ripponstallard} showed that every component of the 
  set $A(f)\subset I(f)$ of ``fast'' escaping points,
  introduced by Bergweiler and Hinkkanen \cite{walteraimo},
  is unbounded. 
  For functions in class $\B$, their ideas can be used to obtain 
  the more precise statement of the following
  theorem. This theorem will be used only in Section
  \ref{sec:sinj}.

 \begin{thm}[Existence of unbounded connected sets]
  \label{thm:ripponstallard}
  Let $f\in\B$, let $F$ be a fundamental domain of $f$, and let $R>0$. 
  Then there exists an unbounded closed connected set $L$ such that,
  for all $j\geq 0$,
   \[ 
      f^j(L) \subset F\cap I(f)\cap \{z\colon |z|\geq R\}. \]
 \end{thm}
 \begin{remark}
  Since the preparation of this article, there have been several improvements
   on this result. For example, it follows from \cite{eremenkoproperty} that
   the set $L$ can be chosen to be forward invariant. Also, in
   \cite{walterphilgwyneth} it is shown that any tract of any entire 
   transcendental function
   (not necessarily in class $\B$) contains an unbounded closed connected
   set of points which escape within this tract. For the reader's convenience,
   we will nonetheless include the simple direct
   proof of the above theorem here.
 \end{remark}
 \begin{proof}  Let $C_F$ be the constant from Lemma 2.3, and let
  $K$ again be as defined at the beginning of the section. 
  We can
  choose
  $C>C_F$ large enough such that
  $|z| > \max(K+1,R)$ whenever $r(z)\geq C$, and such that
  $F$ contains some point 
  $w_0$ with $r(w_0)<C$. We define a sequence
  $(U_j)$ by letting
  $U_0$ be the unbounded connected component of 
   $\{z \in F: r(z) > C\}$, and denoting the
   unbounded connected component of $f(U_j)\cap F$ by 
  $U_{j+1}$ for each $j\geq 0$. 
  (Note that such a component exists and is unique: in fact, 
   by induction the set $U_j$ contains all points of $F$ of sufficiently
   large modulus.)

 Since
  $U_{j+1}\subset f(U_j)$ (and $f|_{U_j}$ is univalent), we can also define
   \[ V_j := (f|_{U_0})^{-1}( (f|_{U_1})^{-1}(\dots
                    (f|_{U_{j-1}})^{-1}(U_j)\dots)) = (f|_F)^{-n}(U_j), \]
  for all $j\geq 0$. Then $V_j$ is connected and
  $f^j:V_j\to U_j$ is a conformal isomorphism.

  By choice of $C$, every point $z\in U_j$ satisfies
  $r(z) > 2^j C$ and $|z| > \max(K+1,R)$. Thus 
  \begin{equation}
   r(f^k(z)) > 2^k\cdot C   \hspace{0.5cm} \text{and}
   \hspace{0.5cm} \label{eqn:growthoforbits}
   |f^k(z)| > R
  \end{equation}
  for every $z\in V_j$ and $k=0,\dots,j$. 

 \begin{claim}
   For every $j\geq 0$, there exists
    $z\in \partial V_j$ with $r(z)=C$. 
 \end{claim}
 \begin{subproof}
  By choice of $C$, the point $w_0\in F$ does not
   belong to $U_j$, and thus there is some point
   $z_0\in F\cap \partial U_j$. Since $|z_0|\geq K+1$, 
   this means that $z_0\in f(F)\cap \partial f(U_{j-1})$. In other words,
   $(f|_F)^{-1}(z_0)\in F\cap \partial U_{j-1}$. 

  Continuing inductively,
   the (unique) point $z\in \partial V_j$ with $f^j(z)=z_0$ satisfies
   $z\in F\cap \partial U_0$, and therefore $r(z)=C$. 
 \end{subproof}

  Since the set $\{z\in\C: r(z)=C\}$ is bounded, 
    $\wt{L} := \bigcap_{j\geq 0} \cl{V_j}$
  is nonempty (as well as invariant, closed and unbounded). By
  (\ref{eqn:growthoforbits}), $\wt{L}$ satisfies all requirements of
  the theorem, except that it need not be connected.
  However, $\wt{L}\cup\{\infty\}$ is compact and connected, and we
  can complete the proof by letting $L$ be any connected component of
  $\wt{L}$.
 \end{proof}


\section{Proof of the main theorem}
  \label{sec:mainpullbacks}

 \begin{defn}[Extendable sets]
  Let $f\colon \C\to\C$ be an entire function.
   $A\subset\C$ is called an 
   \emph{extendable set} (for $f$) if, whenever
   $z$ tends to $\infty$ in $A$, $|f(z)|$ also tends to $\infty$.
   (In particular, every bounded set is extendable.)
 \end{defn}

 Our main result below will apply in any case where a Siegel disk,
  fixed ray etc.\ is an extendable set. 
  In order
  to obtain the theorems as stated in the introduction,
  we remark the following. 

\begin{figure} 
\subfigure[$f(z)=z\exp(z)$]{\fbox{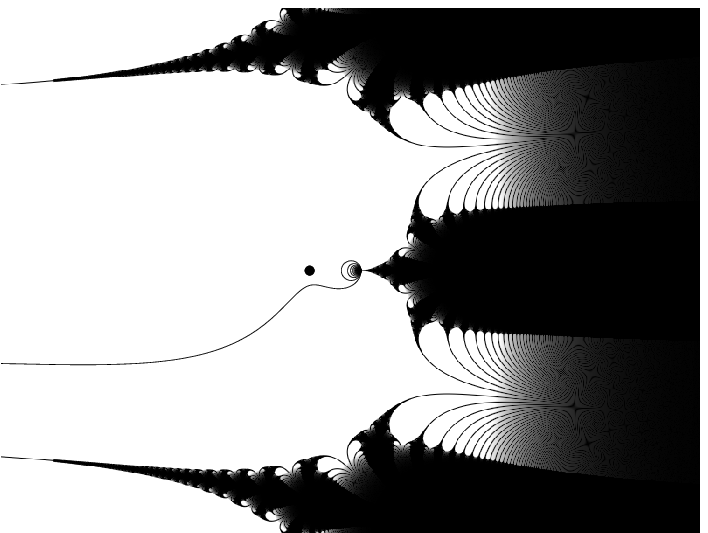}%
  \hspace{3.2pt}%
              \fbox{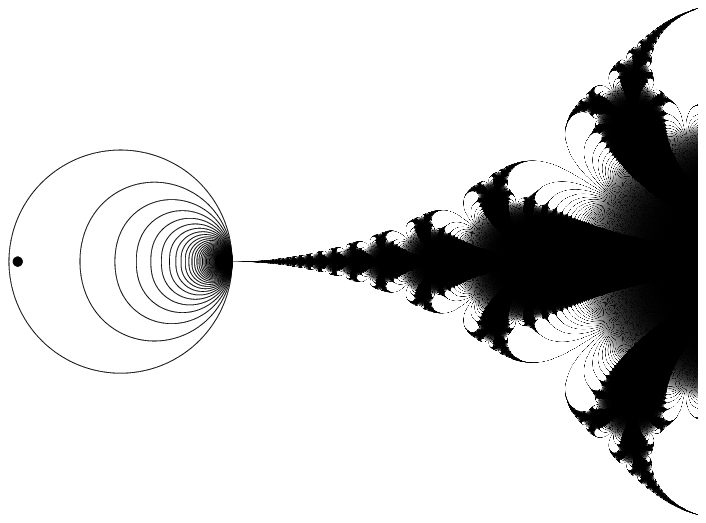}\label{fig:zexpz}}
\subfigure[$f(z)=1/4((z+1)\exp(z)-1)$]{%
   \fbox{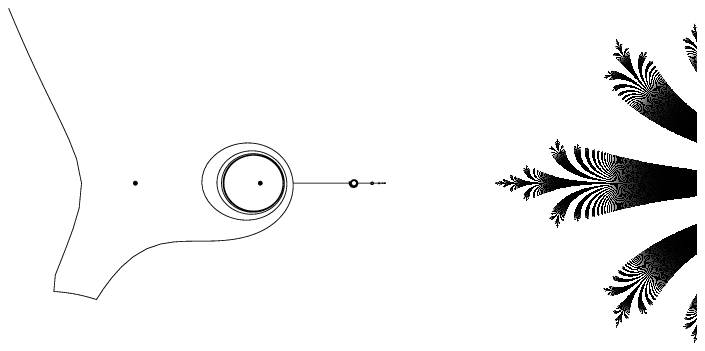}\hspace{3.5pt}%
              \fbox{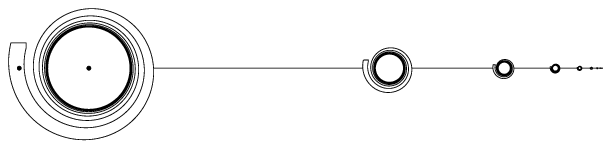}\label{fig:bijective}}
 \caption{Examples of closed connected invariant
   sets $A'$ which are mapped bijectively
   but not homeomorphically
   by an entire function $f$, as described in
   Remark 2 after Observation \ref{obs:extendability}. In both cases,
   the compact
   image $f(A')\subset A'$ is shown in greater magnification on the right.}
\end{figure}

 \begin{obs}[Sufficient conditions for extendability]
  \label{obs:extendability}
  Let $f\colon \C\to\C$ be entire and let $A\subset\C$.
   \begin{enumerate}
    \item $A$ is extendable if and only if $\cl{A}$ is.
       \label{item:closure}
    \item If $f|_A$ is injective, then $A$ is extendable if
      and only if $\partial A$ is.
       \label{item:boundary}
    \item If $A$ is closed and $f\colon A\to\cl{f(A)}$ is a proper map 
     (in particular, if $f\colon A\to\cl{f(A)}$ is a homeomorphism),
     then $A$ is extendable.
       \label{item:homeo} 
   \end{enumerate}
 \end{obs}
 \begin{proof} (\ref{item:closure}) is trivial, and
    (\ref{item:homeo}) follows from the fact that every
    proper map between two topological spaces
    extends continuously to a map between their
    one-point compactifications.

   To prove (\ref{item:boundary}), suppose that
    $\partial A$ is extendable and let
    $(z_n)$ be a sequence in $A$ such that 
    $f(z_n)\to w_0\in \C$. We need to show that
    $z_n$ is bounded. By hypothesis, there exists some $R>0$ such that
    $|f(z)|>|w_0|+1$ for all $z\in\partial A$ with $|z|\geq R$. Let $J$
    be some Jordan curve in $\C\setminus f^{-1}(w_0)$ which surrounds
    $\D_R(0)$. If all but finitely many
    $z_n$ are surrounded by $J$, there is nothing
    to prove. 
 
   Otherwise, set $\eps := \min\bigl(1,\dist(w_0,f(J))\bigr)$ and let $z_k$ be 
    a point on the outside of $J$ with $f(z_k)\in \D_{\eps}(w_0)$. 
    Then the component $U$ of $f^{-1}(\D_{\eps}(w_0))$ containing
    $z_k$ does not intersect $J$, and is thus contained in
    $\C\setminus \D_R(0)$. By choice of $R$, it follows that 
    $U\cap \partial A=\emptyset$, and thus $U\subset A$. Since
    $f|_U$ is injective, it follows easily that
    $f\colon U\to \D_{\eps}(w_0)$ is a conformal map, and thus
    $z_n\to (f|_U)^{-1}(w_0)$ by injectivity of $f|_A$. This means
    that $(z_n)$ is bounded, as required. \end{proof}
 \begin{remark}[Remark 1]
  If $A$ is closed and connected, $A$ contains no critical points 
   and $f\colon A\to\cl{f(A)}$ is 
   a homeomorphism, then it is easy to see that
   there exists a $\C$-neighborhood $U$ of
   $A$ such that $f|_U$ is univalent. For this reason, we
   call a set $A$ satisfying the above assumptions 
   a \emph{set of univalence}. 
   If $f\in\B$ and $A$ is a set of univalence for $f$, then 
   one can show that
   the branch
   $\phi := (f|_{A})^{-1}$ can actually be defined on a large domain of 
   a particularly nice form, namely one whose complement is a 
   union of finitely many arcs to $\infty$ and finitely many
   compact connected sets. However, we do not require this fact
   in this article. 
 \end{remark}
 \begin{remark}[Remark 2]
  Suppose again that 
   $A$ is closed, connected and contains no critical points. 
   If $f$ is a polynomial, or if $A$ is bounded, then $A$ is 
   a set of univalence if and only $f|_A$ is injective. This
   is, of course, no longer true for entire functions: consider
   e.g.\ 
   $f := \exp$ and $A:= \R$. For an example where
   $f\colon A\to\cl{f(A)}$ is bijective but not a homeomorphism, 
   consider the map $f(z) := z\exp(z)$. 
   Here $0$ is both a parabolic
   fixed point and an asymptotic value. Let
   $B$ be a Jordan curve through $0$ which surrounds the critical
   value $-1/e$. If we let $A$ be the curve obtained
   by analytic continuation of
   the branch $\phi$ of $f^{-1}$ with $\phi(0)=0$ along $B$, then
   $A$ is a Jordan arc from $0$ to $\infty$ and
   $f\colon A\to B$ is bijective. If $B$ was chosen to be e.g.\ 
   the boundary of an
   attracting petal at $0$, then we can set
   $A' := \bigcup_{j\geq 0} f^n(A)$. This set $A'$ is closed and
   forward invariant, and 
   $f\colon A'\to \cl{f(A')}$ is bijective
   (see Figure \ref{fig:zexpz}).

  This last example is perhaps not quite satisfactory since the set
   $A'$ contains a singular value. This can be avoided
   by a simple modification. Indeed, let
    $f(z) := 1/4\bigl((z+1)e^z - 1\bigr)$ (which is obtained from our
    previous example by affine coordinate changes in the domain and
    range).  Here the critical point
    $c=-2$ and the asymptotic value $a=-1/4$ both belong to the
    immediate basin of $0$ (which is the entire Fatou set). We can
    let $B$ consist of a small circle around $a$, together with a
    curve spiralling in towards this circle in both directions, and
    surrounding the critical value $f(c)$ (see Figure
    \ref{fig:bijective}). Then, as above, 
    there is an unbounded set $A$ which is mapped bijectively to
    $B$, and if $B$ was chosen correctly, then $A$ and $B$ are disjoint.
    Connecting $A$ and $B$ by an interval of the real axis, and
    adding all forward iterates as well as the fixed point $0$, 
    we obtain a closed invariant set
    $A'$ which contains no singular values and for which 
    $f\colon A'\to \cl{f(A')}$ is bijective;
    see Figure \ref{fig:bijective}.
 \end{remark}

 We are now ready to prove Theorem
  \ref{mainthm:pullbacks}. In fact, we show the following more
  general result. 

 \begin{thm}[Invariant extendable sets]
  \label{thm:mainthmpullbacks}
  Let $f\in \B$ and let $A\subset\C$ be an extendable
   set. 
 \begin{enumerate}
  \item Suppose that $A$ is connected and
   $f(A)\subset A$. Then for every $R>0$, 
   there exists some
    $R'>0$ such that 
    \[ A\cap\{|z|\geq R'\} \subset
         \{z\in I(f)\colon |f^n(z)|\geq R \text{ for all $n\geq 0$} \}. \]
    \label{item:invariantmaintheoremstatement}
  \item More generally, suppose that
    $C\subset A\cap f^{-1}(A)$ is connected. 
    Let $R>0$, and let $X$ denote the set of all
    $z\in C$ with the following property:
    \begin{center} 
      If $n\geq 0$ with $|f^m(z)|\geq R$ for $m=0,\dots,n$,
       then $f^n(z)\in C$.
    \end{center}

   Then there exists $R'>0$ such that
    \[ X\cap\{|z|\geq R'\} \subset
        \{z\in I(f)\colon
           |f^n(z)|\geq R \text{ for all $n\geq 0$} \}. \]
   \label{item:complicatedmaintheoremstatement}
  \end{enumerate}
 \end{thm}
 \begin{proof} 
  (\ref{item:invariantmaintheoremstatement}) follows from
   (\ref{item:complicatedmaintheoremstatement}) by setting
   $C := A$. 
   Thus, it suffices to prove
   (\ref{item:complicatedmaintheoremstatement}). 
   We claim that there is $T>0$ such that
   $\{z\in C\colon |z|\geq T\}$ is contained in finitely many 
   fundamental domains of $f$. 

  To prove this,
   set $K := 1+\max\bigl(|f(0)|,\max_{s\in S(f)}|s|\bigr)$ as before, 
   and let $\gamma$ be the curve to $\infty$
    used in the definition of fundamental domains. Recall that
    $f(\gamma)\subset \cl{\D_K(0)}$. 
    $A$ is extendable, so we can find
    $T\geq T_1>0$ such that, for all $z\in A$,
     \begin{align*}
       |z|\geq T_1 &\Longrightarrow |f(z)|>K \hspace{0.5cm}\text{and}\\
       |z|\geq T &\Longrightarrow |f(z)|\geq T_1.
     \end{align*}
    If $z\in A\cap f^{-1}(A)$ with $|z|\geq T$, then
    $|f(z)|>K$, and thus $z$ belongs to some tract of
    $f$. However, also $|f(f(z))|>K$, and in particular
    $f(z)\notin\gamma$; thus $z$ belongs to some fundamental domain
    of $f$. In other words, every component of $C\setminus D_T(0)$ is contained
    in a single fundamental domain.   
     We may suppose that $T$ was chosen large enough such that
    $C\cap \D_{T}(0)\neq\emptyset$. Since
    $C$ is connected, it follows that
    every fundamental domain which meets
    $C\setminus \D_{T}(0)$ intersects the circle $\{|z|=T\}$.
    As remarked in Section 
    \ref{sec:expansion}, 
    there are only finitely many fundamental domains with
    this property.

  Let $F_1,\dots, F_k$ denote these fundamental domains.  
   Let $z\in X$ and let $n\geq 0$ with
   $|f^m(z)|\geq \max(T,R)$ for all $m\in\{0,\dots,n\}$. Then, 
   by definition
   of $X$, 
   \[ f^n(z)\in C\cap \{|z|\geq T\} \subset \bigcup F_j. \] 
  By Lemma \ref{lem:pointsgrowinfundamentaldomains}, 
   there exists $R'>0$ such that
  \[ X\cap \{|z| \geq R'\} \subset
     \{z\in I(f)\colon |f^n(z)|\geq \max(T,R)\text{ for all $n$}\}. 
       \qedhere \]
 \end{proof}

Since Siegel disks never contain escaping points, we have the following immediate corollary. 

 \begin{cor}[Extendable Siegel disks are bounded]
  \label{cor:univalentsiegel}
  Let $f\in \B$, let $U$ be a Siegel disk of $f$ and suppose that
   $U$ is an extendable set. Then 
   $U$ is bounded.  \qedoutsideproof
 \end{cor}

\begin{proof}[Proof of Theorem \ref{mainthm:univalentsiegel}]
  If $f:\partial U\to\partial U$ is a homeomorphism, this implies that
   $\cl{U}$ (and hence $U$) is an extendable set by Observation
   \ref{obs:extendability}. 
   Theorem \ref{mainthm:univalentsiegel}
    now follows from Corollary \ref{cor:univalentsiegel}.
\end{proof}

 \begin{cor}[Bounded accumulation sets]
   \label{cor:univalentrays}
  Let $f\in \B$ and let $\gamma\colon(-\infty,1]\to\C$ 
   be a curve with $f(\gamma(t))=\gamma(t+1)$ for all $t\leq 0$. If 
   $\gamma$ is an extendable set, then $\gamma$ is bounded.
 \end{cor}
 \begin{proof} We will apply Theorem
   \ref{thm:mainthmpullbacks} (\ref{item:complicatedmaintheoremstatement})
   to the sets $A:=\gamma\bigl((-\infty,1]\bigr)$ and 
   $C:= \gamma\bigl((-\infty,0]\bigr)$. 
   If $z\in C$ such that
   $f(z)\notin C$, then $f(z)\in \gamma\bigl([0,1]\bigr)$.
   Thus, if we set $R := 1 + \max_{t\in[0,1]}|\gamma(t)|$, then
   the set $X$ as defined in Theorem
   \ref{thm:mainthmpullbacks}  
   (\ref{item:complicatedmaintheoremstatement}) is all of $C$. 
   Since every point of
   $C$ eventually maps to $\gamma\bigl([0,1]\bigr)$ under iteration, 
   this means that there is $R'>0$ such that 
      $C \subset \cl{\D_{R'}(0)}$. \end{proof}

\section{Obtaining Extendability}
  \label{sec:sinj}

 It is often difficult to show that a given set is extendable. 
 In this
  section, we describe how to do this in the case where 
  all singular values (with the possible exception of one critical value)
  are contained in the Julia set. The main idea is to
  apply Theorem \ref{thm:ripponstallard}
  and the following
  observation. 

 \begin{lem}[Obtaining extendability]
   \label{lem:extendability}
  Let $f\in\B$. Let $A\subset\C$ be a connected
   set and set $B:= \cl{f(A)}$. 
  \begin{enumerate}
   \item Suppose that there exists some connected open set $G\subset A$
     such that $f|_G$ is univalent and 
     $B\cap S(f)\subset f(G)$. Suppose furthermore that
     $\infty$ is accessible from every component $W$ of 
     $\C\setminus B$ for which $W\cap S(f) \neq \emptyset$. 
     Then $f\colon \cl{A}\to B$ is a homeomorphism. 
     \label{item:allaccessible}
   \item Suppose that
     there exists some isolated point
     $s_0\in S(f)$ which is a critical (but not an asymptotic) value such that
     $B\cap (S(f)\setminus\{s_0\})=\emptyset$ and such that
     $\infty$ is accessible from every component $W$ of 
     $\C\setminus B$ with $W\cap (S(f)\setminus \{s_0\}) \neq \emptyset$. 
     Then $A$ is an extendable set.
     \label{item:onecriticalvalue}
  \end{enumerate}
 \end{lem}
 \begin{proof} Set $S_0 := f(G)$ for 
   (\ref{item:allaccessible}) and $S_0 := \{s_0\}$ for
   (\ref{item:onecriticalvalue}). 
   Let
   $\W$ denote the set of those components $W$ of 
   $\C\setminus B$ which satisfy 
   $W\cap (S(f)\setminus S_0) \neq \emptyset$. 
  Note that $\mathcal{W}$ is an open cover of the compact set
      $S(f)\setminus S_0$ and contains no proper subcover. Hence
      $\mathcal{W}$ is finite. 

  Now, for each component $W\in\W$, we can choose some compact, connected 
   and full (i.e., non-separating)
   subset $K_W\subset W$ with
   $S(f)\cap W\subset K_W$, 
   and a curve $\gamma_W\subset W\setminus K_W$ with one endpoint
   in $K_W$ and the other at $\infty$. 
   Then the set 
   $V := \C\setminus \bigcup_{W\in\W} (K_W\cup\gamma_W)$ is 
   simply connected, contains $B$ and is disjoint from
   $S(f)\setminus S_0$. Let $U$ be the component of
   $f^{-1}(V)$ which contains $A$. 

  In case of (\ref{item:allaccessible}), we claim that
   $f\colon U\to V$ is a conformal isomorphism. Indeed, if
   $G=\emptyset$ (which is the only case in which we will apply
   this lemma), then this follows by the monodromy theorem.
   Otherwise, it is not difficult to see that
   the branch $(f|_G)^{-1}$ extends to a branch of $f^{-1}$ on
   $V$, whose image is necessarily $U$.

  In case of (\ref{item:onecriticalvalue}), 
   $f\colon U\to V$ is either conformal or a finite-degree covering 
   with a single branched point. In either case,
   $f\colon \cl{A}\to B$ is a proper map, which implies that
   $\cl{A}$ (and thus $A$) is extendable. \end{proof}

 \begin{thm}[Accessibility of infinity]
   \label{thm:sinj}
  Let $f\in\B$ and let $U\subset\C$ be connected. Suppose
   that there are three
   fundamental domains $F_1,F_2,F_3$ of $f$ such that 
   \[ U\cap \{z\in I(f)\colon f^k(z)\in F_i
                  \text{ for all
               sufficiently large $k$}\}=\emptyset \]
   for all $i\in\{1,2,3\}$.
   Then $\infty$ is accessible from every component $W$ of
    $\C\setminus \cl{U}$ which satisfies $W\cap J(f)\neq \emptyset$. 
 \end{thm}
 \begin{proof} By Theorem
   \ref{thm:ripponstallard}, for each $i\in\{1,2,3\}$,  
   there is an unbounded closed connected set
   $C_i\subset I(f)$ with
   $f^k(C_i)\subset F_i$ for all $k\geq 0$.

  Let $z_0\in W\cap J(f)$, and set $\delta := \dist(z_0,\cl{U})$ and
   $D := \D_{\frac{\delta}{2}}(z_0)$. 
   Let
   $i\in \{1,2,3\}$. Then, since $z_0\in J(f)$,
   there exists some large $n$ such that
   $f^n(D)\cap C_i\neq \emptyset$. Let
   $A_i'$ be a component of
   $f^{-n}(C_i)$ with $A_i'\cap D\neq\emptyset$. 

  Note that the three sets
   $A_i'$
   are unbounded and pairwise disjoint. For each $i$,
   let $A_i$ be the closure of
   some unbounded component of
   $A_i'\setminus D$. (Such a component exists by a simple
   application of the
   boundary bumping theorem
   \cite[Theorem 5.6]{continuumtheory}.) These components are disjoint from
   $U$ by the assumption on $U$. 

  It follows that there exists $j\in\{1,2,3\}$ with
   $A_j \subset W$ (see Lemma \ref{lem:threesets}). By
   Lemma \ref{lem:infinityaccessible}, $\infty$ is accessible from $W$. \end{proof}

 \begin{proof}[Proof of Theorem \ref{mainthm:nonsingularsiegel}]
  If $U$ is a Siegel disk, then $U$ does not intersect
   $I(f)$ and therefore satisfies the
   assumption of Theorem \ref{thm:sinj}. Thus $\infty$ is
   accessible from every component of $\C\setminus\cl{U}$
   which intersects the Julia set. If $S(f)\subset J(f)$ and
   $S(f)\cap\partial U = \emptyset$, then we can apply
   Lemma \ref{lem:extendability} (\ref{item:allaccessible})
   to see that 
   $f\colon \cl{U}\to \cl{U}$ is a homeomorphism. 
 \end{proof}

 \begin{proof}[Proof of Theorem \ref{mainthm:nonsingularcosinesiegel}]
  Let
   $U_0\mapsto U_1\mapsto\dots\mapsto U_m\mapsto U_0$ be a cycle of
   Siegel disks for a map $f$ with two critical values and no asymptotic
   values.
   Suppose that, for each $j$, at least one
   of the two critical values of $f$ belongs to 
   $J(f)\setminus \partial U_j$.
   (This assumption is automatically satisfied if no $\partial U_j$ contains
    a critical value. Indeed, each Siegel disk boundary is contained in 
    the postcritical set and
    $f$ has no wandering domains \cite{alexmisha}, so
    at least one critical value must belong to the Julia set.) 

   Applying Theorem \ref{thm:sinj} and Lemma \ref{lem:extendability}
   (\ref{item:onecriticalvalue}), we see that $U_j$ is an extendable
   set for each $j$. Thus $U_0$ is an extendable set for
   $f^m$, and 
   the claim follows from Corollary \ref{cor:univalentsiegel}. 
 \end{proof}

 For cases in which it may not be possible to control the
  eventual behavior of points, let us also show the
  following variant of Theorem \ref{thm:sinj}. 

 \begin{thm}[Accessibility when $S(f)\subset J(f)$]
  \label{thm:sinj2}
  Let $f\in\B$ and suppose that $S(f)\subset J(f)$. 
  Then for every $R>0$ and $\eps>0$, there exists an $n_0$ 
  with the following property. If
   \[ U\subset \{z\in\C\colon \exists n\geq n_0\colon |f^n(z)| <  R\} \]
  is connected with
   $\dist(U,S(f))\geq \eps$, then $\infty$ is accessible
   from every component of $\C\setminus \cl{U}$. 
 \end{thm}
 \begin{proof} Applying Theorem
   \ref{thm:ripponstallard} to three different fundamental domains of $f$,
   there exist three unbounded closed connected sets
   \[
     C_1,C_2,C_3\subset 
     \{z\in I(f)\colon |f^m(z)|\geq R \text{ for all $m$}\} \]
   such that $f^m(C_i)\cap f^n(C_j) = \emptyset$ whenever
    $i\neq j$ and $m,n\geq 0$.
   Let 
  \[ \K := \bigcup_{s\in S(f)} \cl{\D_{\frac{\eps}{2}}(s)}. \]
  As in the proof of Lemma \ref{lem:extendability}, this set
   has finitely many components. For each
   such component $K$ and each
   $i\in\{1,2,3\}$, we can again find some
   $n_{K,i}$ and a closed unbounded set
   $A_{K,i}$ connecting $K$ to $\infty$ such that
   $f^{n_{K,i}}(A_{K,i})\subset C_i$.

  Let $\displaystyle{n_0 := \max_{K,i} n_{K,i}}$. Then every set $U$ as in the
   statement of the theorem is disjoint from 
   $\K\cup \bigcup_{K,i} A_{K,i}$,
   and the proof proceeds as in Theorem \ref{thm:sinj}. \end{proof}

 \begin{cor}[Nonsingular rays are bounded] \label{cor:sinjrays}
  Let $f\in\B$ with $S(f)\subset J(f)$. Suppose that 
   $\gamma\colon (-\infty,1]\to\C$ 
   is a curve with $f(\gamma(t))=\gamma(t+1)$ for all $t\leq 0$ and
   $\cl{\gamma}\cap S(f)=\emptyset$. Then there exists $t\leq 0$ such that 
   \[ f\colon \cl{\gamma\bigl((-\infty,t]\bigr)}\to
        \cl{\gamma\bigl((-\infty,t+1]\bigr)} \]
   is a homeomorphism. 
   In particular, $\gamma$ is bounded.   
 \end{cor}
 \begin{proof} Let $R:=1+\max_{t\in[0,1]}|\gamma(t)|$ and
   $\eps := \dist(\gamma,S(f))$. By Theorem \ref{thm:sinj2}, there
   exists $t := -n_0$ such that
   $\gamma\bigl((-\infty,t]\bigr)$
   satisfies the assumptions of Lemma
   \ref{lem:extendability} (\ref{item:allaccessible}), 
   proving that $f$ is a homeomorphism on its closure.  
   That $\gamma$ is bounded now follows from 
   Corollary \ref{cor:univalentrays} and 
   Observation \ref{obs:extendability}. \end{proof}

\section{Landing of Periodic Rays} \label{sec:landing}

 The classical snail lemma
  \cite[Lemma 16.2]{jackdynamicsthird} states that the landing point of
  an invariant curve for a holomorphic mapping $f$ cannot be an
  irrationally indifferent fixed point. We will prove a generalization
  of this fact which allows us to prove the landing of certain
  periodic rays. Our methods are quite similar to those used by Perez-Marco
  in his study of hedgehogs \cite{perezmarco}; since first
  preparing this article
  we have also learned that Risler \cite{rislerinvariant} also studied completely 
  invariant compact sets of univalent functions using
  similar considerations.

 Let us say
  that a pair $(f,K)$ of a compact set $K$ and a holomorphic map
  $f$ defined in a neighborhood $U$ of $K$ 
  \emph{has the snail lemma property} if every curve
  $\gamma\colon (-\infty,1]\to U\setminus K$ with
  $f(\gamma(t))=\gamma(t+1)$ and 
  $\lim_{t\to-\infty}\dist(\gamma(t),K)=0$ lands at a 
  repelling or parabolic fixed point of $f$.

 \begin{lem}[Univalent snail lemma for the circle]
  \label{lem:snailunivalentcircle}
  Let $U$ be a neighborhood of $S^1$, and let
   $f\colon U\to\C$ be a univalent function with $f(S^1)=S^1$.
   Then $(f,S^1)$ has the snail lemma property.
 \end{lem}
 \begin{proof} Let $\gamma$ be a curve as in the 
   snail lemma property. Then 
   $f$ is not of finite order, as otherwise no such curve
   $\gamma$ can exist. We may also assume (by reflection in $S^1$)
   that $\gamma\subset\D$ and (by restriction of $U$) that
   $f$ has no fixed points outside the unit circle. 
   Let us consider two cases. 

 \smallskip
 \noindent 
  \emph{First case: $f$ possesses at least one fixed point}. 
  Since $f\neq \id$, the number of fixed points of $f$ is finite.
  The unit circle
  $S^1$ is invariant under $f$, and thus these fixed points
  must be either attracting, repelling or parabolic. They
  cut the circle into finitely many intervals, and 
  points in such an interval converge to one
  endpoint under forward iteration
  and to the other under backwards iteration. Thus, every
  interval, with the exception of one endpoint, is contained in the
  basin of attraction (or repulsion) of the other endpoint. 

  The accumulation set of $\gamma$ cannot intersect any of the
   basins of attraction, because every point on $\gamma$ eventually maps 
   to $\gamma\bigl([0,1]\bigr)$. Since the accumulation set is connected,
   it consists of a single repelling or parabolic point, as required.

\medskip
 \noindent
  \emph{Second case: $f$ has no fixed points.}
   The argument in this case is completely
   analogous to the proof of the 
   classical snail lemma. We shall therefore omit some
   of the details in the proof. 

  We may assume that $0\notin U$.
   Set $\wt{U}:= \exp^{-1}(U\cap\D)$ and let $\wt{\gamma}\subset \wt{U}$ 
   be any lift of
   $\gamma$ under $\exp$. We can then choose a lift
   $\wt{f}$ of $f$ such that $\wt{f}(\wt{\gamma}(t)) = 
    \wt{\gamma}(t+1)$. 

  Since $f$, and thus $\wt{f}$, has no fixed points,
   \[ c := \max_{r\in\R} |\im \wt{f}(ir) - r| > 0. \]
  Choose $\eps>0$ such that $|\im \wt{f}(z) - \im z|\geq \frac{c}{2}$ whenever
   $|\re z| \leq \eps$.    
  To fix ideas, let us assume that $\im \wt{f}(ir) < r$ for all
   $r\in\R$. It follows
   that $\theta(t)\to +\infty$ as $t\to-\infty$, where
    $\theta(t) := \im\wt{\gamma}(t)$. 

      We can thus pick some $t_0<0$ such that 
       $|\re \wt{\gamma}(t)| \leq \eps$ for $t\leq t_0$ and 
       $\theta(t)\leq \theta(t_0)$ for $t\geq t_0$. 
   Let 
   $\wt{V}$ be the component of 
    \[ \{z\in\wt{U}\setminus\wt{\gamma}: \im z > \theta(t_0)\} \]
   whose boundary contains the line $\{ir\colon r\geq \theta(t_0)\}$. 
   It follows easily that $\wt{V}\subset \wt{f}(\wt{V})$.

  Therefore $V := \exp(\wt{V})$ is a one-sided neighborhood of $S^1$ with
   $f^{-1}(V)\subset V$.
   Since the boundary of $V$ is contained 
   in $U$, the iterates $f^{-n}|_V$ converge locally uniformly
   to a fixed point of $f$ by
   \cite[Lemma 5.5]{jackdynamicsthird}. This contradicts our assumption.    
  \end{proof}

 \begin{lem}[General univalent snail lemma]
  \label{lem:snailunivalent}
  Let $K\subset\C$ be compact and connected, and suppose that
   $f$ is a function univalent in a neighborhood $U$ of $K$, with
   $f(K)=K$. Then $(f,K)$ has the snail lemma property.
 \end{lem}
 \begin{proof} Let 
   $\gamma\subset U\setminus K$
   be a curve as in the definition of the snail lemma property and let
   $V$ be the component of $\C\setminus K$ containing $\gamma$. Then
   $U\cap V$ is connected and is thus mapped by 
   $f$ into some component of $\C\setminus K$. 
   Since
   $f(U\cap V)$ contains $\gamma$ and thus intersects $V$, it follows that
   $f(U\cap V)\subset V$.

  Let $\phi\colon V\to\D$ be a Riemann mapping of $V$, and define
  \[ g\colon\phi(U\cap V)\to \D; z\mapsto \phi(f(\phi^{-1}(z))). \]
  Since $f$ is continuous in a neighborhood of $K$ and
   $f^{-1}(K)=K$, every prime end of $K$ is mapped to a prime
   end of $K$ by $f$. Thus $g$ extends continuously to $S^1$ by
   Carath\'eodory's Theorem \cite[Theorem 2.15]{pommerenke}.
   By the Schwarz Reflection Principle
   \cite{ahlforscomplexanalysis}, 
   $g$ extends to an analytic function on a neighborhood
   of $S^1$.

  This extended function $g$ is univalent, and
   thus the curve
   $\phi(\gamma)$ lands at a repelling or parabolic fixed point $z_0$ of
   $g$ by Lemma \ref{lem:snailunivalentcircle}. 
   Let $D$ be a linearizing neighborhood or repelling
   petal of $z_0$ which is compactly contained in the 
   domain of definition of $g$ and which
   contains some end piece of $\phi(\gamma)$.

  Then
   $\phi^{-1}(D)$ is invariant under $f^{-1}$ and contains an end 
   piece of $\gamma$. Again, the functions
   $(f|_D)^{-n}$ converge locally uniformly to a
   fixed point of $f$ in $K$ by
   \cite[Lemma 5.5]{jackdynamicsthird}. It follows that $\gamma$ lands
   at this fixed point, which is repelling or parabolic
   by the classical Snail Lemma. 
  \end{proof}

 \begin{cor}[Landing of univalent rays]
   \label{cor:raysgeneral}
  Let $f$ be an entire function and suppose
   that $\gamma\colon (-\infty,1]\to I(f)$ is a curve with
   $f(\gamma(t))=\gamma(t+1)$. Suppose furthermore
   that $\gamma$ does not accumulate at any critical point of
   $f$ and that, for some $t\leq 0$, the restriction
    \[ f \colon \cl{\gamma\bigl((-\infty,t]\bigr)} \to
           \cl{\gamma\bigl((-\infty,t+1]\bigr)} \]
   is a homeomorphism.
   Then $\gamma$ lands at a repelling or parabolic
   fixed point of $f$.
 \end{cor}
 \begin{proof} By Corollary \ref{cor:univalentrays}, 
   the accumulation set $K$ of $\gamma$
   is bounded, and therefore contains no escaping points.
   In particular, $\gamma\cap K=\emptyset$. Since
   $f\colon K\to K$ is a homeomorphism and since $K$ contains no critical points
   of $f$, it follows easily that there exists some neighborhood
   $U$ of $K$ such that $f|_U$ is univalent (see 
   \cite[Lemma 3]{hermansiegel}). The claim now follows from
   Lemma \ref{lem:snailunivalent}. 
   \end{proof}

 \begin{proof}[Proof of Theorem \ref{mainthm:landing}]   
  If we are in the first case of Theorem
   \ref{mainthm:landing}, then Corollary \ref{cor:univalentrays}
   implies (using Observation \ref{obs:extendability}) that
   the accumulation set of $\gamma$ is bounded. In the second case,
   the same follows from Corollary 
   \ref{cor:sinjrays}. So in either case Corollary \ref{cor:raysgeneral}
   implies that $\gamma$ lands at a repelling or parabolic fixed point of $f$,
   as claimed.  
 \end{proof}

  In \cite{landing2new},
   it was shown that periodic rays of exponential 
   maps \emph{always} land, but this requires the ``lambda-lemma'' and 
   deep results on exponential parameter
   space. We can now deduce a special case of
   this theorem without requiring any parameter-space
   arguments, but also without any
   a priori assumptions on hyperbolic expansion. 
   (Compare Appendix \ref{app:landing1}).

 \begin{cor}[Nonsingular exponential rays]
  Suppose that $f(z)=\exp(z)+\kappa$ and that
   $\gamma\colon (-\infty,\infty)\to\C$ is a periodic dynamic ray with
   $\kappa\notin\bigcup_j\cl{f^j(\gamma)}$. 
   Then $\gamma$ lands at a repelling or
   parabolic periodic point of $f$. 
 \end{cor}
 \begin{proof} First suppose that $\kappa\in F(f)$. Since
   $f$ (like all maps with only finitely many singular values)  
   has no wandering domains
   \cite{alexmisha}, this implies that $\kappa$ belongs to the
   basin of an attracting or parabolic
   periodic point. In this case, all periodic rays of $f$ land
   by Corollary \ref{cor:landing1}.

  So now suppose that $\kappa\in J( f)$. Let $n$ be the period of $\gamma$.
   Let $N$ be large enough and set 
  \[ g_j := \cl{f^j\bigl(\gamma\bigl((-\infty,-N]\bigr)\bigr)}
    \]
  By Theorem \ref{thm:sinj2} and Lemma
   \ref{lem:extendability}, if $N$ was chosen large enough,
   then $f\colon g_j\to g_{j+1}$ is a homeomorphism for $j=0,\dots,n$.
   Therefore
   $f^n\colon g_0\to g_n$ is a homeomorphism, and
   the claim follows by applying Corollary
   \ref{cor:raysgeneral} to $f^n$. \end{proof}   

\appendix

\section{Two topological facts} \label{app:topology}
 This section is dedicated to proving the two simple topological facts
  which were
  used in Section \ref{sec:sinj}.

 \begin{applem}[Accessibility criterion] \label{lem:infinityaccessible}
  Let $U$ be a domain in $\C$. Suppose that there exists an unbounded
  closed connected set $A\subset U$. Then $\infty$ is accessible from $U$.
 \end{applem}
\begin{proof}
  For each $z\in A$, let $\delta(z):=\min(\dist(z,\partial U),1)$. For
  every $n\in\N$, the set
    \[ A_n := \{z\in A\colon n-1\leq |z|\leq n\} \]
  is compact; thus there exists a finite set $K_n\subset A_n$ such
  that
    \[ A_n \subset \bigcup_{z\in K_n} \D_{\delta(z)}(z). \]
  We claim that there is an infinite sequence 
   $z_1,z_2,z_3,\dots\in K := \bigcup_j K_j$ with the property that
  $\D_{\delta(z_j)}(z_j)\cap\D_{\delta(z_{j+1})}(z_{j+1})\neq\emptyset$ for all
  $j$ and such that
   $z_j\neq z_{j'}$ for all $j\neq j'$. 
   Indeed, consider the (infinite)
   graph on the set $K$ in which $z$ and $w$ are adjacent if
   $\D_{\delta(z)}(z)\cap \D_{\delta(w)}(w)\neq\emptyset$. Then the
   graph $G$ is connected and contains arbitrarily long paths. By
   K\"onig's lemma (see e.g.~\cite[Lemma 7.1.3]{diestel}), $G$ contains
   an infinite path, as
   desired.

  So let $z_1,z_2,\dots$ be a sequence as above. Then the curve obtained
  by connecting every $z_j$ to $z_{j+1}$ by a straight line segment 
  is a curve to $\infty$ in $U$. \end{proof}

 \begin{applem}[Separation lemma]
   \label{lem:threesets}
   Let $D\subset\C$ be a closed disk, and let
     $U\subset \C$ be connected with
     $\dist(U,D)>0$.  

  Suppose that
     $C_1,C_2,C_3$ are pairwise
     disjoint closed connected sets such that
     $C_j\cap U=\emptyset$, $C_j\cap D\neq\emptyset$ and 
     $C_j\setminus D$ is connected for every $j$.
   Then there exists $i\in \{1,2,3\}$ such that
    $C_i\cap \cl{U}=\emptyset$.
 \end{applem}
 \begin{proof} 
  If, for some $i_1,i_2\in \{1,2,3\}$, the sets
   $U$ and $C_{i_1}\setminus D$ belong to different
   components of $\C\setminus (C_{i_2}\cup D)$, then 
   $C_{i_1} \cap \cl{U} = \emptyset$,
   and we are done.
   Otherwise, let $K_i := \Ch\setminus G_i$, where $G_i$
   is the component of $\C\setminus (C_i
            \cup D)$ containing $U$. 
   Then each $K_i$ is a nonseparating continuum, and any two of
   these
   intersect exactly in $D\cup\{\infty\}$. 
   
  The set
   $\Ch\setminus \bigcup K_i$ has exactly three components 
   $W_1, W_2,W_3$ \cite{straszewicz}, which we may suppose labelled such that
   $U\subset W_1$. It is a simple application of 
   Janiczewski's theorem (compare \cite[Theorem 1.9]{pommerenkeunivalent})
   that there exist 
   $i,j\in\{1,2,3\}$, $i\neq j$, such that
   $K_i\cup K_j$ does not separate $W_2$ and $W_3$. By relabelling,
   we may assume that
   $i=1$ and $j=2$. 

  Let $V$ be the component of 
   $\Ch\setminus (K_1\cup K_2)$ which contains $W_2$ and $W_3$, 
   and define
   $K_3' := K_3\setminus (D\cup\{\infty\})$.
   By construction, $K_3'$ is connected,   intersects $V$
   and is
   disjoint from
   $K_1\cup K_2$. Thus $K_3\subset V$.
   Also note that 
   $V\cap W_1 = \emptyset$, since $K_1\cup K_2$ separates the
   Riemann sphere by \cite{straszewicz}.
   Therefore
   \[ C_3 \subset K_3'\cup D \subset V \cup D \subset
        (\C\setminus \cl{W_1}) \cup D \subset
           \C\setminus \cl{U}. \qedhere \]  
  \end{proof}

\section{Hyperbolic Contraction} \label{app:landing1}
  For completeness, let us discuss here in which situations hyperbolic
   contraction can be used to show the landing of fixed rays of
   an entire function. The proof in the polynomial case 
    (Corollary \ref{cor:polylanding})
   was given by
   Douady and Hubbard \cite{orsay}. The following is something of a
   ``folk theorem''; special cases with essentially the same proof
   can be found
   e.g.\ in \cite{nuriastandard,thesis,expattracting,expper}. 
  \begin{appthm}[Landing of fixed rays via contraction]
   \label{thm:landing1}
     Let $X$ be a Riemann surface, and let $U\subset X$ be a hyperbolic
      domain. Let $h:U\to X$ be holomorphic such that $h(U)\supset U$ and
      $h:U\to h(U)$ is a covering map. Furthermore suppose that $h$ is
      not an irrational rotation of a disk, punctured disk or annulus,
      and that $h^n\neq \id$ for $n\geq 1$.  

    Let $\gamma\colon (-\infty,1]\to U$ be a
     curve
     with $h(\gamma(t))=\gamma(t+1)$
     for all $t\geq 1$. 
   Then every accumulation point of
     $\gamma$ in 
     \[ U\cup \bigl\{z\in\partial U\colon h \text{ extends continously to
                        $U\cup\{z\}$} \bigr\}  \]
     is a fixed
    point of (the extension of) $h$.
%
  \end{appthm}
  \begin{proof} 
  Denote hyperbolic distance and length in $U$ 
  (compare \cite[Chapter 2]{jackdynamicsthird}) by
  $\dist_U$ and $\ell_U$, respectively.
  Since $h$ is a covering map, $h$ is a local isometry 
  between $U$ and $h(U)$
  (with their respective hyperbolic metrics). 
  By the Schwarz Lemma, the inclusion $U\to h(U)$ does not expand
  the hyperbolic metric, and so $h$ is either a local isometry (if $U=h(U)$)
  or otherwise strictly expands the hyperbolic
  metric of $U$.

 Let us first suppose that 
  $\gamma|_{[0,1]}$ is smooth; at the end of the proof
  we will sketch how to deal with the general 
  case. 
  For $t\leq 0$, let us denote by $\ell_t$ the hyperbolic
  length of $\gamma([t,t+1])$ in $U$. Since pullbacks under $h$
  contract the hyperbolic metric, 
  $\ell_t$ is a nonincreasing function of $t$.

  Suppose that $z_0\in\partial U$ is an accumulation point of
  $\gamma$; say $\gamma(t_n)\to z_0$ with $t_n\to-\infty$. Then the
  hyperbolic distance between $z_n := \gamma(t_n)$ and
  $h(z_n)=\gamma(t_n+1)$ is at most $\ell_{t_n}\leq \ell_0$; since the
  hyperbolic metric of $U$ blows up near $\partial U$, 
  it follows that the distance between $z_n$ and $h(z_n)$ in $X$ 
  tends to $0$. If $h$ can be
  extended continuously into $z_0$, then $|z_0-h(z_0)|=0$
  by continuity. 

  Now suppose that $h(U)\supsetneq U$ and let 
   $z_0$ be an accumulation point of $\gamma$ 
   in $U$. Again, choose $t_n\to-\infty$ such that 
     $z_n := \gamma(t_n)\to z_0$; we may suppose that
   $t_{n+1}\leq t_n-1$ for all $n$.
  Let $D$ be the closed hyperbolic disk 
   \[ D:=\{z\in U\colon \dist_U(z,z_0)\leq 2\ell_0\}. \]
  Since $h$ strictly expands the hyperbolic metric, 
  there exists a number $\lambda > 1$ such
  that
  $\|\operatorname{D}\hspace{-0.5mm}h(z) \|_{\operatorname{hyp}}\geq \lambda$
  for all $z\in U$ with $h(z)\in D$; in fact, we can take 
  $\lambda := \max_{z\in D} \frac{\rho_U(z)}{\rho_{h(U)}(z)}$
  (where $\rho_U$ and $\rho_{h(U)}$ are the densities of the corresponding
   hyperbolic metrics). 

  If $n_0$ is large enough, then
   $\dist_U(\gamma(t_n),z_0)<\ell_0$ for $n\geq n_0$, and thus
   $\gamma\bigl([t_n,t_n+1]\bigr) \subset D$. Therefore 
    \begin{align*}
       \ell_{t_n} =
       \ell_U\bigl(\gamma[t_n,t_n+1]\bigr) &=
       \ell_U\bigl( h(\gamma[t_n-1,t_n])\bigr) \\ &\geq
       \lambda \ell_U\bigl(\gamma([t_n-1,t_n])\bigr) =
       \lambda \ell_{t_n-1} \geq \lambda \ell_{t_{n+1}}. 
    \end{align*}
  By induction, $\ell_{t_n} \leq \lambda^{n_0-n}\ell_{t_{n_0}}$ for
   $n\geq n_0$, which means that $\ell_t\to 0$ as $t\to-\infty$. 
    Therefore
   \[ \dist_U(z,h(z))=\lim \dist_U(z_n,h(z_n)) = \lim \ell_{t_n} = 0,
   \]
  as required. 

  Finally, suppose that $h(U)=U$. Then by assumption and by 
   \cite[Theorem 5.2]{jackdynamicsthird}, the iterates of $f$ tend
   locally 
   uniformly to $\infty$ (in the sense that they leave every compact subset
   of $U$).
, or we are in one of the exceptional cases from the statement of
   our theorem. In the former case, it follows easily that $\gamma$ has
   no accumulation points in $U$ at all. 

  If $\gamma|_{[0,1]}$ is not smooth, then the proof can easily be modified
   in a number of ways. For example, let us define $d(t_1,t_2)$ to be the
   (hyperbolic) length of the unique hyperbolic geodesic of $U$ homotopic
   to $\gamma|_{[t_1,t_2]}$. Then we can set 
     \[ \ell(t) := \max_{s\leq t} d(s,s+1) = 
         \max_{s\in [t-1,t]} d(s,s+1), \] 
   and the proof goes through as before. (Alternatively, replace
   $\gamma|_{[0,1]}$ by a smooth curve $\tilde{\gamma}$
   in the same homotopy class
   and continue $\tilde{\gamma}$ using pullbacks. Using hyperbolic contraction
   arguments as above, it is easy to see that
   both curves will have the same accumulation sets.)
      \end{proof}

  \begin{appthm}[Fixed rays do not land at infinity]
    \label{thm:nolandingatinfty}
   Let $f\in\B$, and suppose that
    $\gamma\colon (-\infty,1]\to\C$ is a curve with
    $f(\gamma(t)) = \gamma(t+1)$ for all $t\leq 0$. 
    If $z_0\in\Ch$ with 
    $\lim_{t\to-\infty} \gamma(t)=z_0$, then 
    $z_0\neq \infty$.
  \end{appthm}
  \begin{proof} Suppose, by contradiction, that $z_0=\infty$. 
    Then $\gamma$ is an extendable set, which is
    impossible by
    Corollary \ref{cor:univalentrays}. \end{proof}

  \begin{appcor}[Landing criterion]
    \label{cor:landingcriterion}
   Let $f\in\B$, and let $\gamma\colon (-\infty,1]\to I(f)$ with
    $f(\gamma(t))=\gamma(t+1)$. Then $\gamma$ lands at a
    repelling or parabolic fixed point of $f$ if and only if
    there exists some domain $U$ such that
    $U\subset  f(U)$, $f\colon U\to f(U)$ is a covering map and
    $\gamma\bigl((-\infty,T]\bigr)\subset U$ for some $T<0$. 
  \end{appcor}
  \begin{proof} The ``if'' part is a combination of
   Theorem \ref{thm:landing1}, Theorem
   \ref{thm:nolandingatinfty} and the classical Snail Lemma.
   The ``only if'' part is trivial since we can let $U$ be
   a linearizing neighborhood or repelling petal. \end{proof}

  \begin{appcor}[Rays not intersecting the postsingular set]
   \label{cor:landing1}
   Let $f\in\B$ and $\gamma$ as above, and suppose that
    $\gamma\bigl((-\infty,T]\bigr)\cap \P(f) = \emptyset$ for some
    $T<0$. Then $\gamma$ lands at a repelling or parabolic
    fixed point of $f$. 
  \end{appcor}
 \begin{proof} Let $V := \C\setminus \P(f)$, and let $U$ be the component of
    $f^{-1}(V)$ containing $\gamma\bigl((-\infty,T-1]\bigr)$. 
    The claim follows from the landing
    criterion. 
  \end{proof}

  \begin{appcor}[Polynomial rays]
   \label{cor:polylanding}
   Let $f$ be a polynomial. Then every periodic ray of $f$ lands 
    at a repelling or parabolic fixed point of $f$. 
  \end{appcor}
  \begin{proof} Let $U$ be the basin of
   infinity and apply Theorem \ref{thm:landing1}. 
 \end{proof}

  Both Corollary \ref{cor:landingcriterion} and
   Corollary \ref{cor:raysgeneral} give necessary and sufficient conditions
   for fixed rays to land. However, in practice it appears to be often
   difficult to obtain these conditions. Let us consider the
   family $\Ek\colon z\mapsto \exp(z)+\kappa$ 
   of exponential maps as an example. Suppose that, for some
   $\Ek$, there exists a fixed ray $\gamma$ which accumulates on the
   singular value $\kappa$, whose orbit in turn accumulates on all of
   $\gamma$. (It can easily happen that the singular orbit is dense in 
   $\C$; in fact this behavior is generic in the bifurcation locus
   \cite[Theorem 5.1.6]{thesis}.) 

  It is easy to see that, in this situation, no set $U$ as in
   Corollary \ref{cor:landingcriterion} could exist. Thus, in order to
   prove landing of $\gamma$ by a hyperbolic contraction result, we need
   to \emph{a priori} show that $\gamma$ does not accumulate at $\kappa$.
   This would appear to be extremely difficult without any prior
   dynamical information: even in well-controlled cases there 
   are dynamic rays which do not land, and 
   whose accumulation sets are actually indecomposable
   continua
   \cite{nonlanding}. 

  In the exponential family, this problem can be circumvented by using
   parameter space arguments. This gives some hope that more general
   landing theorems are also true in higher-dimensional parameter spaces
   such as the space of cosine maps, but 
   unfortunately the methods in the exponential
   case break down completely. 
   On the other hand, Corollary \ref{cor:raysgeneral} shows that
   failure of univalence really is
   the only possible obstruction to landing of
   periodic rays.

\nocite{perezmarco}
\nocite{hermansiegel}
\nocite{rippon}
\bibliographystyle{hamsplain}
\bibliography{M:/Latex/Biblio/biblio}

\end{document}